\documentclass[11pt]{jdg-p-1-1}

\newtheorem{theorem}{Theorem}[section]
\newtheorem{lemma}[theorem]{Lemma}

\newtheorem{proposition}[theorem]{Proposition}

\theoremstyle{definition}

\newtheorem{example}[theorem]{Example}

\theoremstyle{remark}
\newtheorem{remark}[theorem]{Remark}

\def \b {\beta}

\def\lf{\left}
\def\ri{\right}
\def\bbar{{\bar{\beta}}}
\def\a{\alpha}

\def\g{\gamma}
\def\e{\epsilon}
\def\p{\partial}
\def\delbar{{\bar{\delta}}}

\def\C{\Bbb C}
\def\R{\Bbb R}

\def\k{\Upsilon}

\def\ba{{\bar{\alpha}}}
\def\bb{{\bar{\beta}}}
\def\abb{\alpha{\bar{\beta}}}

\def\Q{{\mathcal Q}}

\def \D {\Delta}

\numberwithin{equation}{section}

\begin{document}

\title{A note on  Perelman's LYH inequality}

\author{Lei Ni}
\address{Department of Mathematics, University of California at San Diego, La Jolla, CA 92093}

\email{lni@math.ucsd.edu}
\thanks{The author was supported in part by NSF Grants  and
an Alfred P. Sloan Fellowship, USA}


\date{September 2005}


\begin{abstract} We give a proof to the Li-Yau-Hamilton type inequality
claimed by Perelman on the fundamental solution to the conjugate
heat equation. The rest of the paper is devoted to improving the
known differential inequalities of Li-Yau-Hamilton type via
monotonicity formulae.
\end{abstract}

\maketitle

\section{Introduction}

In  \cite{P}, Perelman {\it proved} a Li-Yau-Hamilton type (also
called differential Harnack) inequality for the fundamental solution
of the conjugate heat equation, in the presence of the Ricci flow.
More precisely, let $(M, g_{ij}(t))$ be a solution to Ricci flow:
\begin{equation}\label{1}
\frac{\p}{\p t} g_{ij}=-2R_{ij}
\end{equation}
on $M\times [0, T]$ and let $H(x,y, t)=\frac{e^{-f}}{(4\pi
\tau)^{\frac{n}{2}}}$ (where $\tau=T-t$) be the fundamental solution
to  the conjugate heat equation $u_\tau-\Delta u+Ru=0$. (More
precisely we should write the fundamental solution as $H(y,t; x,
T)$, which satisfies $ \left(-\frac{\p}{\p t}+\D_{y}+R(y,
t)\right)H=0$ for any $(x, t)$ with $t<T$ and $\lim_{t\to T} \int_M
H(y, t; x, T)f(y, t)\, d\mu_t(y)=f(x, T)$.) Define
$$
v_H=\left[\tau \left(2\Delta f-|\nabla f|^2+R\right)+f-n\right]H.
$$
Here all the differentiations are  taken with respect to $y$, and
$n=\dim_\R(M)$. Then $v_H\le 0$ on $M\times [0, T]$. This result is
a differential inequality of Li-Yau type \cite{LY}, which
 has important consequences in the later part of \cite{P}.  For example it is essential
  in proving the pseudo-locality theorems.
 It is also
 crucial in localizing the entropy formula  \cite{N3}.

In Section 9 of \cite{P},  the following important differential
equation
\begin{equation}\label{2}
\left(\frac{\p}{\p \tau}-\D
+R\right)v_u=-2\tau|R_{ij}+\nabla_i\nabla_j
f-\frac{1}{2\tau}g_{ij}|^2u \end{equation} is stated for any
positive solution $u$ to the conjugate heat equation, whose
integration on $M$ gives the celebrated entropy formula for the
Ricci flow. One can consult various resources (e.g. \cite{N1}) for
the detailed computations of this equation, which can also be done
through a straightforward calculation, after knowing the result.
\cite{P} then proceeds the proof of the claim $v_H\le 0$ in a clever
way by checking that for any $\tau_*$ with $T\ge \tau_*>0$, $\int_M
v_H(y) h(y)\, d\mu_{\tau_*}(y)\le 0$, for any smooth function
$h(y)\ge 0$ with compact support. In order to achieve this, in
\cite{P} the heat equation $\left(\frac{\p}{\p t}-\D\right)h(y,
t)=0$ with the `initial data' $h(y, T-\tau_*)=h(y)$ (more precisely
$t=T-\tau_*$), the given compactly supported nonnegative function,
is solved. Applying (\ref{2}) to $u(y, \tau)=H(x, y, \tau)$, one can
easily derive as in \cite{P}, via integration by parts, that
\begin{equation}\label{3}
\frac{d}{ d\tau}\int_M v_H h\, d\mu_\tau=-2\int_M
\tau|R_{ij}+\nabla_i\nabla_j f-\frac{1}{2\tau}g_{ij}|^2Hh\,
d\mu_\tau\le 0.
\end{equation}
The Li-Yau type inequality  $v_H\le 0$  then follows from the above
monotonicity, provided the claim that \begin{equation}\label{4}
\lim_{\tau \to 0}\int_M v_Hh\, d\mu_\tau \le 0. \end{equation}

The main purpose of this note is to prove (\ref{4}), hence provide a
complete proof of the claim $v_H\le 0$.  This will be done in
Section 3 after some preparations in Section 2. It was written in
\cite{P} that `{\it it is easy to see}' that $\lim_{\tau\to 0}\int_M
v_Hh\, d\mu_\tau = 0$. It turns out that the proof found here need
to use some gradient estimates for positive solutions, a quite
precise estimate on the `reduced distance' and the monotonicity
formula (\ref{3}). (We shall focus on the proof of (\ref{4}) for the
case when $M$ is compact and leave the more technical details of
generalizing it to the noncompact setting to the later refinements.)
Indeed the claim that $\lim_{\tau\to 0}\int_M v_Hh\, d\mu_\tau = 0$
follows from a blow-up argument of \cite{P}, after we have
established (\ref{4}). Since our argument is a bit involved, this
may not be {\it the proof}.

In Section 4 we derive several monotonicity formulae, which improve
various Li-Yau-Hamilton inequalities for linear heat equation
(systems) as well as  for Ricci flow, including the original
Li-Yau's inequality. In Section 5 we illustrate  the localization of
them by applying a general scheme of \cite{EKNT}.

\section{Estimates and results needed} We shall collect some known
results and derive some estimates needed for proving (\ref{4}) in
this section. We  need the asymptotic behavior of the fundamental
solution to the conjugate heat equation for small $\tau$. Let
$d_\tau(x,y)$ be the distance function with respect to the metric
$g(\tau)$. Let  $B_\tau(x, r)$ ($Vol_\tau$) be the ball of radius
$r$ centered at $x$ (the volume) with respect to the metric
$g(\tau)$.

\begin{theorem} \label{Theorem 1} Let $H(x,y, \tau)$ be the fundamental solution to the ({\it
backward in $t$}) conjugate heat equation. Then as $\tau \to 0$ we
have that
\begin{equation}\label{5}
H(x,y,\tau)\sim \frac{\exp\left(-\frac{d_0^2(x,y)}{4\tau}\right)}
{\left(4\pi\tau\right)^{\frac{n}{2}}}\sum_{j=0}^{\infty}\tau^ju_j(x,y,\tau).
\end{equation} By (\ref{5}) we mean that exists $T>0$ and sequence $u_j\in
C^{\infty}(M\times M\times [0, T])$ such that
$$
H(x,y,\tau)-\frac{\exp \left(-\frac{d_0^2(x,y)}{4\tau}\right)}
{\left(4\pi\tau\right)^{\frac{n}{2}}}\sum_{j=0}^{k}\tau^ju_j(x,y,\tau)=w_k(x,y,
\tau)
$$
with
$$
w_k(x,y, \tau)=O\left(\tau^{k+1-\frac{n}{2}}\right)
$$
as $\tau \to 0$, uniformly for all $x, \, y\in M$.  The function
$u_0(x, y, \tau)$ can be chosen so that $u_0(x, x, 0)=1$.
\end{theorem}

This result was proved in details, for example in \cite{GL}, when
there is no zero order term $R(y, \tau) u(y, \tau)$ in the equation
$\frac{\p}{\p \tau}u-\D u+ Ru=0$ and replacing $d_0(x,y)$ by
$d_\tau(x,y)$. However, one can check that the argument carries over
to this case if one assumes that the metric $g(\tau)$ is $C^\infty$
near $\tau=0$. One can consult \cite{SY, CLN} for intrinsic
presentations.

Let
$$
\mathcal {W}_h(g, H, \tau)=\int_M v_Hh\, d\mu_\tau
$$
where $h$ is the previously described solution to the heat equation.
It is clear that for any $\tau$ with $T\ge \tau>0$,
$\mathcal{W}_h(g, H, \tau)$ is a well-defined quantity. A priori it
may blow up as $\tau\to 0$. It turns out that in our course of
proving that $\lim_{\tau\to 0} \mathcal{W}_h(g, H,\tau)\le 0$ we
need to  show first that exists $C>0$, which may depends on the
geometry of the Ricci flow solution $(M, g(\tau))$ defined on
$M\times [0, T]$, but independent of $\tau$ (as $\tau\to 0$) so that
$\mathcal{W}_h(g, H, \tau)\le C$ for all $T \ge\tau>0$. The
following lemma supplies the key estimates for this purpose.

\begin{lemma}\label{Lemma 2} Let $(M, g(t))$ be a smooth solution to the Ricci
flow on $M\times[0, T]$. Assume that there exist $k_1\ge 0$ and
$k_2\ge0$, such that  the Ricci curvature $R_{ij}(g(\tau))\ge -k_1
g_{ij}(\tau)$ and $\max(R(y, \tau), |\nabla R|^2(y, \tau))\le k_2$,
on $M\times [0, t]$.

(i) If $u\le A$ is a positive solution to the conjugate heat
equation on $M\times [0, T]$, then there exists $C_1$ and $C_2$
depending on $k_1$, $k_2$ and $n$ such that for $0< \tau\le \min(1,
T)$,
\begin{equation}\label{6}
\tau\frac{|\nabla u|^2}{u^2}\le \left(1+C_1\tau\right)\left(\log
\left(\frac{A}{u}\right)+C_2\tau\right)\end{equation}

(ii) If  $u$ is a positive solution to the conjugate heat equation
on $M\times[0, T]$, then there exists $B$, depending on $(M,
g(\tau))$ so that for $0\le\tau\le \min(T, 1)$,
\begin{equation}\label{7}
\tau \frac{|\nabla u|^2}{u^2}\le \left(2+C_1\tau\right)\left(\log
\left(\frac{B}{u\tau^{\frac{n}{2}}}\int_M u\,
d\mu_{\tau}\right)+C_2\tau\right). \end{equation}
\end{lemma}
\begin{remark}\label{Remark1} Here and thereafter we use the same $C_i$ ($B$) at the
different lines  if they just differ only by a constant depending on
$n$. Notice that $\int_M u\, d\mu_{\tau}$ is independent od $\tau$
and equals to $1$ if $u$ is the fundamental solution. The proof to
the lemma given below is a modification of some  arguments in
\cite{H}.
\end{remark}
\begin{proof} Direct computation, under a unitary frame, gives
\begin{eqnarray*}
 \left(\frac{\p}{\p \tau}-\Delta \right)\left(\frac{|\nabla
u|^2}{u}\right) &=&
-\frac{2}{u}\left|u_{ij}-\frac{u_iu_j}{u}\right|^2+\frac{|\nabla
u|^2}{u}
R\\
&\, &+\frac{-4R_{ij}u_i u_j -2\langle \nabla (R u), \nabla
u\rangle}{u}\\
&\le& (4+n)k_1\frac{|\nabla u|^2}{u}+2|\nabla R||\nabla u|\\
&\le& \left[(4+n)k_1+1\right]\frac{|\nabla u|^2}{u}+k_2u
\end{eqnarray*}
and
\begin{eqnarray*}
\left(\frac{\p}{\p \tau}-\Delta
\right)\left(u\log\left(\frac{A}{u}\right)\right)&=&\frac{|\nabla
u|^2}{u}+Ru-Ru\log\left(\frac{A}{u}\right)\\
&\ge &\frac{|\nabla u|^2}{u}-nk_1u-k_2u\log
\left(\frac{A}{u}\right).
\end{eqnarray*}
Combining the above two equations together we have that
$$
\left(\frac{\p}{\p \tau}-\Delta \right)\Phi\le 0
$$
where $$\Phi=\varphi\frac{|\nabla
u|^2}{u}-e^{k_2\tau}u\log\left(\frac{A}{u}\right)-(k_2+nk_1e^{k_2})\tau
u$$ with $\varphi=\frac{\tau}{1+\left[(4+n)k_1+1\right]\tau}$, which
satisfies
$$
\frac{d}{d\tau} \varphi +\left[(4+n)k_1+1\right]\varphi <1.
$$
By the maximum principle we have that
$$
\varphi\frac{|\nabla u|^2}{u}\le
e^{k_2\tau}u\log\left(\frac{A}{u}\right) +(k_2+nk_1e^{k_2})\tau u.
$$
From this one can derive (\ref{6}) easily.

To prove the second part, we claim that for $u$, a positive solution
to the conjugate heat equation, there exists a $C$ depending on $(M,
g(\tau))$ such that \begin{equation}\label{claim}  u(y, \tau)\le
\frac{C}{\tau^{\frac{n}{2}}}\int_M u(z, \tau)\, d\mu_\tau(z).
\end{equation}
This is a mean-value type inequality, which can be proved via,  for
example the Moser iteration. Here we follow \cite{H}. We may assume
that $\sup_{y\in M, 0\le \tau\le 1}\tau^{\frac{n}{2}}u(y,\tau)$ is
finite. Otherwise we may replacing $\tau$ by $\tau_\e =\tau-\e$ and
let $\e\to 0$ after establishing the claim for $\tau_\e$. Now let
$(x_0, \tau_0)\in M\times [0,1]$ be such a space-time point that
$\max \tau^{\frac{n}{2}}u(y, \tau) =\tau_0^{\frac{n}{2}}u(y_0,
\tau_0)$. Then we have that
$$
\sup_{M\times [\frac{\tau_0}{2}, \tau_0]}u(y, t)\le
\left(\frac{2}{\tau_0}\right)^{\frac{n}{2}}\tau_0^{\frac{n}{2}}u(y_0,
\tau_0)=2^{\frac{n}{2}}u(y_0, t_0).
$$
Noticing this upper bound, we  apply (\ref{6}) to $u$ on
$M\times[\frac{\tau_0}{2}, \tau_0]$, and conclude  that
$$
\frac{\tau_0}{2}\left(\frac{|\nabla u|^2}{u^2}\right)(y, \tau_0)\le
(1+C_1\tau_0)\left(\log\left( \frac{2^{\frac{n}{2}}u(y_0,
\tau_0)}{u(y, \tau_0)}\right)+C_2\tau_0\right).
$$
Let $g=\log\left( \frac{2^{\frac{n}{2}}u(y_0, \tau_0)}{u(y,
\tau_0)}\right)+C_2\tau_0$.  The above can be written as
$$
|\nabla \sqrt{g}|\le \sqrt{\frac{1+C_1\tau_0}{2\tau_0}}
$$
which implies that
$$
\sup_{B_{\tau_0}\left(y_0,
\sqrt{\frac{\tau_0}{1+C_1\tau_0}}\right)}\sqrt{g}(y, \tau_0)\le
\sqrt{g}(y_0, \tau_0)+\frac{1}{\sqrt{2}}.
$$
Rewriting the above in terms of $u$ we have that
$$
u(y, \tau_0)\ge 2^{\frac{n}{2}}u(y_0,
\tau_0)e^{-\left(\frac{1}{2}+\frac{2}{\sqrt{2}}\sqrt{\frac{n}{2}\log
2+C_2}\right)}=C_3 u(y_0, \tau_0)
$$
for all $y\in B_{\tau_0}\left(y_0,
\sqrt{\frac{\tau_0}{1+C_1\tau_0}}\right)$. Here we have also used
$\tau_0\le 1$. Noticing that
$$
Vol_{\tau_0}\left( B_{\tau_0}\left(y_0,
\sqrt{\frac{\tau_0}{1+C_1\tau_0}}\right)\right)\ge C_4
\tau_0^{\frac{n}{2}}
$$
for some $C_4$ depending on the geometry of $(M, g(\tau_0))$.
Therefore we have that
$$
\frac{C_5}{\tau_0^{\frac{n}{2}}}\int_M u(y, \tau_0)\,
d\mu_{\tau_0}(y)\ge u(y_0, \tau_0)
$$
for some $C_5$ depending on $C_3$ and $C_4$. By the way we choose
$(y_0, \tau_0)$ we have that
$$
 \tau^{\frac{n}{2}}u(y, \tau)\le \tau_0^{\frac{n}{2}}u(y_0,
\tau_0) \le C_5\int_M u(y, \tau_0)\, d\mu_{\tau_0}(y) = C_5\int_M
u(y, \tau)\, d\mu_{\tau}(y).
$$
This proves the claim (\ref{claim}). Now the estimate (\ref{7})
follows from (\ref{6}), applying to $u$ on $M\times[\frac{\tau}{2},
\tau]$, and the just proved (\ref{claim}), which ensures the needed
upper bound for applying the estimate (\ref{6}).
\end{proof}

If $u=\frac{e^{-f}}{\left(4\pi\tau\right)^{\frac{n}{2}}}$ is the
fundamental solution to the conjugate heat equation we have that
$\int_M u\, d\mu_\tau=1$. Therefore, by (\ref{7}), we have that
\begin{equation}\label{8}
\int_M \tau|\nabla f|^2 uh \, d\mu_\tau \le
\left(2+C_1\tau\right)\int_M \left(\log B +f +C_2\tau\right)uh\,
d\mu_\tau. \end{equation} On the other hand, integration by parts
can rewrite
\begin{eqnarray*}
\mathcal{W}_h(g, u, \tau)&=&\int_M \tau|\nabla f|^2uh\,
d\mu_\tau-2\tau\int_M \langle \nabla f, \nabla h\rangle u\,
d\mu_\tau \\
&\quad& +\tau \int_M R uh\, d\mu_\tau  +\int_M (f-n)uh\,
d\mu_\tau\\&=&I+II+III+IV.
\end{eqnarray*}
The $I$ term can be estimated by (\ref{8}), whose right hand side
contains only one `bad' term $\int_M fuh\, d\mu_\tau$ in the sense
that it could possibly blow up. The second term
$$
II=2\tau\int_M \langle \nabla u, \nabla h\rangle \,
d\mu=-2\tau\int_M u\D h \, d\mu_\tau
$$
 is clearly bounded as $\tau \to 0$. In fact $II\to 0$ as $\tau\to
0$. The same conclusion
 obviously holds for $III$. Summarizing above,  we reduce the
question of bounding from above the quantity $\mathcal{W}_h(u, g,
\tau)$ to bounding one single term
$$
V=\int_M fuh\, d\mu_\tau
$$
from above (as $\tau \to 0$). We shall show later that $
\lim_{\tau\to 0} V\le 0.$ To do this we need to use the `reduced
distance', introduced by Perelman in \cite{P} for the Ricci flow
geometry.

Let $x$ be a fixed point in $M$. Let  $\ell(y, \tau)$ be the reduced
distance in \cite{P}, with respect to $(x, 0)$ (more precisely
$\tau=0$). We collect the  relevant  properties of $\ell(y, \tau)$
in the following lemma (Cf. \cite{Ye, CLN}).

\begin{lemma}\label{Lemma 3} Let $\bar L(y, \tau)=4\tau \ell(y, \tau)$.

(i) Assume that there exists a constant $k_1$ such that
$R_{ij}(g(\tau))\ge -k_1g_{ij}(\tau)$, $\bar L(y, \tau)$ is a local
Lipschitz function on $M\times[0, T]$;

(ii) Assume that there exist constant $k_1$ and $k_2$ so that
$-k_1g_{ij}(\tau)\le R_{ij}(g(\tau))\le k_2g_{ij}(\tau)$. Then
\begin{equation}\label{9}
\bar L(y, \tau)\le e^{2k_2\tau}d^2_0(x, y)+\frac{4k_2 n}{3}\tau^2
\end{equation} and
\begin{equation}\label{10}
d^2_0(x, y)\le e^{2k_1\tau}\left(\bar L(y, \tau)+\frac{4k_1
n}{3}\tau^2\right); \end{equation}

(iii) \begin{equation}\label{11}\left(\frac{\p}{\p \tau}-\D
+R\right)\left(\frac{\exp\left(-\frac{\bar L(y,
\tau)}{4\tau}\right)}{(4\pi\tau)^{\frac{n}{2}}}\right)\le 0.
\end{equation}
\end{lemma}

\begin{proof} The first two claims follow from the definition by straight forward
checking. For (iii), it was proved in Section 7 of  \cite{P}. By now
there are various resources where the detailed proof can be found.
See for example \cite{Ye} and \cite{CLN}.
\end{proof}

The consequence of (\ref{9}) and (\ref{10}) is that
$$
\lim_{\tau \to 0}\frac{\exp\left(-\frac{\bar L(y,
\tau)}{4\tau}\right)}{(4\pi\tau)^{\frac{n}{2}}}=\delta_x(y),
$$
which together with (\ref{11}) implies that $H$, the fundamental
solution to the conjugate heat equation, is bounded from below as
$$H(x, y, \tau)\ge \frac{\exp\left(-\frac{\bar L(y,
\tau)}{4\tau}\right)}{(4\pi\tau)^{\frac{n}{2}}}.$$ Hence
\begin{equation}\label{12}
f(y, \tau)\le \frac{\bar{L}(y, \tau)}{4\tau}. \end{equation} This
was proved in \cite{P} out of the claim $v_H\le 0$. Since we are in
the middle of proving $v_H\le 0$, we have to show the above
alternative of obtaining (\ref{12}).

 \section{ Synthesis }

Now we assembly the results in the previous section to prove
(\ref{4}). As the  first step we show that $\mathcal{W}_h(g, H,
\tau)$ is bounded (thanks to the monotonicity (\ref{3}), it is
sufficient to bound it from above) as $\tau \to 0$, where $H(x, y,
\tau)$ is the fundamental solution to the conjugate heat equation
with $H(x, y,0)=\delta_x(y)$. By the reduction done in the previous
section we only need to show that
$$
V=\int_M fuh\, d\mu_\tau
$$
is bounded from above as $\tau \to 0$. By (\ref{12}) we have that
\begin{eqnarray*}
 &\, & \limsup_{\tau \to 0}\int_M fHh\, d\mu_\tau \le \limsup_{\tau \to
0}\int_M \frac{\bar L (y, \tau)}{4\tau} H(x, y, \tau) h(y, \tau)\, d\mu_\tau(y)\\
&\quad &\quad  \le \limsup_{\tau \to 0}\int_M \frac{d_0^2(x,
y)}{4\tau}H(x, y, \tau)
h(y, \tau)\, d\mu_\tau(y)\\
&\quad &\quad +\lim_{\tau \to 0}\int_M
\left(\frac{e^{k_2\tau}-1}{4\tau}d_0^2(x,
y)+\frac{k_2n}{3}\tau\right)H(x, y, \tau) h(y, \tau)\, d\mu_\tau(y).
\end{eqnarray*}
Here we have used (\ref{9}) in the last inequality. By Theorem
\ref{Theorem 1}, some elementary computations give  that
$$
\lim_{\tau \to 0}\int_M \frac{d_0^2(x, y)}{4\tau}H(x, y, \tau) h(y,
\tau)\, d\mu_\tau(y)=\frac{n}{2} h(x, 0).
$$
Since $\frac{e^{k_2\tau}-1}{4\tau}d_0^2(x, y)+\frac{k_2n}{3}\tau$ is
a bounded continuous function even at $\tau=0$, we have that
$$
\lim_{\tau \to 0}\int_M \left(\frac{e^{k_2\tau}-1}{4\tau}d_0^2(x,
y)+\frac{k_2n}{3}\tau\right)H(x, y, \tau) h(y, \tau)\,
d\mu_\tau(y)=0.
$$
This completes our proof on finiteness of $\limsup_{\tau\to 0}
\int_M f H h\, d\mu_\tau.$ In fact we have proved that
\begin{equation}\label{13} \limsup_{\tau\to 0} \int_M (f-\frac{n}{2}) H h\,
d\mu_\tau\le 0. \end{equation} By the just proved finiteness of
$\mathcal{W}_h(g, H, \tau)$ as $\tau\to 0$, and the (entropy)
monotonicity (\ref{3}), we know that the limit $\lim_{\tau \to 0}
\mathcal{W}_h(g, H, \tau)$ exists. Let
$$
\lim_{\tau \to 0} \mathcal{W}_h(g, H, \tau)=\lim_{\tau \to 0} \int_M
v_H h\, d\mu_\tau =\a
$$
for some finite $\alpha$.  Hence  $\lim_{\tau \to 0}
\left(\mathcal{W}_h(g, H, \tau)-\mathcal{W}_h(g, H,
\frac{\tau}{2})\right)=0$. By (\ref{3}) and the mean-value theorem
we can find  $\tau_k\to 0$ such that
$$
\lim_{\tau_k \to 0} \tau_k^2\int_M |R_{ij}+\nabla_i\nabla_j
f-\frac{1}{2\tau_k}g_{ij}|^2Hh\, d\mu_{\tau_k}=0.
$$
By the Cauchy-Schwartz inequality and the H\"older inequality we
have that
$$
\lim_{\tau_k \to 0} \tau_k\int_M \left(R+\D f
-\frac{n}{2\tau_k}\right)Hh\, d\mu_{\tau_k}=0.
$$
This implies that
$$
\lim_{\tau \to 0}\mathcal{W}_h(g, H, \tau)=\lim_{\tau_k\to 0}\int_M
\left(\tau_k(\Delta f -|\nabla f|^2)+f-\frac{n}{2}\right)H h\,
d\mu_{\tau_k}.
$$
Again the integration by parts shows that
\begin{eqnarray*}
\int_M \tau_k(\Delta f -|\nabla f|^2)Hh\, d\mu_{\tau_k}&=&\int_M
\tau_k \langle \nabla H, \nabla h\rangle\,
d\mu_{\tau_k}\\
&=&-\tau_k\int_M H\Delta h\, d\mu_{\tau_k}\to 0.
\end{eqnarray*}
Hence by (\ref{13}) $$ \lim_{\tau \to 0}\mathcal{W}_h(g, H,
\tau)=\lim_{\tau_k\to 0}\int_M (f-\frac{n}{2})Hh\, d\mu_{\tau_k}\le
0.
$$
This proves $\alpha \le 0$, namely (\ref{4}).

The claim that $\alpha=\lim_{\tau \to 0}\mathcal{W}_h(g, H, \tau)=0$
can now be proved by the blow-up argument as in Section 4 of
\cite{P}. Assume that $\alpha<0$. One can easily check that this
would imply that $\lim_{\tau \to 0}\mu(g, \tau)<0$. Here $\mu(g,
\tau)$ is the invariant defined in Section 4 of \cite{P}. In fact,
noticing that $h(y, \tau)>0$ for all $\tau\le \tau_*$ (where the
$\tau_*$ is the one we fixed in the introduction). Therefore by
multiple $\frac{1}{h(x, 0)}$ (more precisely $\frac{1}{h(x, \cdot)}$
at $\tau=0$) to the original $h(y, \tau)$, we may assume that $h(x,
0)=\int_M H(x, y, \tau)h(y, \tau)\, d\mu_\tau=1$. Let $\tilde u(y,
\tau)=H(x,y,\tau)h(y, \tau)$ and $\tilde f=-\log \tilde u
-\frac{n}{2}\log(4\pi)$.  Now direct computation yields that
$$
\mathcal{W}_h( g, H, \tau)=\mathcal{W}(g, \tilde u,  \tau)+\int_M
\left(\tau\left(\frac{|\nabla h|^2}{ h}\right)-h\log h\right)H\,
d\mu_\tau.
$$
Noticing that the second integration goes to $0$ as $\tau\to 0$, we
can deduce that $\mathcal{W}(g, \tilde u,  \tau)<0$ for sufficient
small $\tau$ if $\alpha <0$. This, together with the fact $\int_M
\tilde u\, d\mu=1$,  implies that $\mu(g, \tau)<0$ for sufficiently
small $\tau$.
  Now  Perelman's blow-up argument in the Section 4 of \cite{P} gives a  contradiction with the
sharp logarithmic Sobolev inequality  on the Euclidean space
\cite{G}. (One can consult, for example \cite{N1, STW}, for more
details of this part.)

\begin{remark}\label{Remark 2} The method of proof here follows a similar idea used
in \cite{N1}, where the asymptotic limit of the entropy as $\tau\to
\infty$ was computed. Note that we have to use properties of the
reduced distance, introduced in Section 7 of \cite{P}, in our proof,
while the similar, but slightly easier, claim that $\lim_{\tau\to
0}\mathcal{W}(g, H, \tau)=0$ appears much earlier in Section 4 of
\cite{P}.
\end{remark}

The proof can be easily modified to give the asymptotic behavior of
the entropy defined in \cite{N1} for the fundamental solution to the
linear heat equation, with respect to a fixed Riemannian metric.
Indeed if we restrict to the class of complete Riemannian manifolds
with non-negative Ricci curvature we have the following estimates.
\begin{proposition}\label{heat-est} For any $\delta>0$, there exists
$C(\delta)$ such that
\begin{equation} \label{grad-log3}
\frac{|\nabla H|^2}{H}(x, y, \tau)\le 2\frac{H(x, y,
\tau)}{\tau}\left(C(\delta)+\frac{d^2(x,y)}{(4-\delta)\tau}\right)
\end{equation}
and
\begin{equation}\label{lapla-log2}
\Delta H (x, y, \tau) +\frac{|\nabla H|^2}{H}(x, y, \tau)\le
2\frac{H(x, y,
\tau)}{\tau}\left(C(\delta)+4\frac{d^2(x,y)}{(4-\delta)\tau}\right).
\end{equation}
\end{proposition}
The previous argument for the Ricci flow case can be transplanted to
show that
$$
\tau(2\Delta f -|\nabla f|^2)+f-n\le 0
$$
where $H(y, \tau; x, 0)=\frac{1}{(4\pi \tau)^{n/2}}e^{-f}$ is the
fundamental solution to the heat operator $\frac{\p}{\p \tau}-\D$.
This gives a rigorous argument for the inequality (1.5) (Theorem
1.2) of \cite{N1}, for both the compact manifolds and complete
manifolds with non-negative Ricci (or Ricci curvature bounded from
below). For the full detailed account please see \cite{CLN}.

\section{Improving Li-Yau-Hamilton estimates via monotonicity
formulae}

The proof of (\ref{4}) indicates a close relation between the
monotonicity
 formulae and the differential inequalities of Li-Yau type. The hinge
is simply Green's second identity. This was discussed very generally
in \cite{EKNT}. Moreover if we chose $h$ in the introduction to be
the fundamental solution to the time dependent heat equation
($\frac{\p}{\p t}-\D$) centered at $(x_0, t_0)$ we can have a better
upper bound on $v_H(x_0, t_0)$ in terms of the a weighted integral
which is non-positive.  In fact, this follows from the
representation formula for the solutions to the non-homogenous
conjugate heat equation. More precisely, since $h(y, t; x_0, t_0)$
 is the fundamental
solution to the heat equation (to make it very clear, $v_H$ is
defined with respect to $H=H(y, t; x, T)$,  the fundamental solution
to the conjugate heat equation centered at $(x, T)$ with $T>t_0$),
we have that
$$
\lim_{t\to t_0}\int_M h(y, t; x_0, t_0) v_H(y, t)\, d\mu_t(y)
=v_H(x_0, t_0).
$$
On the other hand from (\ref{2}) we have that (by  Green's second
identity)
$$
\frac{d}{d t} \int_M h v_H\, d\mu_t =2\tau\int_M
|R_{ij}+f_{ij}-\frac{1}{2\tau}|^2H h\, d\mu_t.
$$
Therefore
$$
\lim_{t\to T}\int_M hv_H\, d\mu_t - v_H(x_0,
t_0)=\int_{t_0}^T2\tau\int_M
|R_{ij}+f_{ij}-\frac{1}{2\tau}g_{ij}|^2H h\, d\mu_t\, dt.
$$
Using the fact that $\lim_{t\to T}v_H=0$ we have that
$$
v_H(x_0, t_0)=-2\int_{t_0}^T (T-t)\int_M
\left|R_{ij}+f_{ij}-\frac{1}{2(T-t)}\right|^2Hh\, d\mu_t \, dt \le
0,
$$
which sharpens the estimate $v_H\le 0$ by providing a non-positive
upper bound. Noticing also the duality $h(y, t; x_0, t_0)=H(x_0,
t_0; y, t)$ for any $t>t_0$ we can express everything in terms of
the fundamental solution to the {\it (backward) conjugate heat
equation}.

Below we show a few new monotonicity formulae, which expand the list
of examples shown in the introduction of \cite{EKNT} on the
monotonicity formulae, and more importantly improve the earlier
established Li-Yau-Hamilton estimates in a similar way as the above.

For the simplicity let us just consider the K\"ahler-Ricci flow case
even though often the discussions are also valid for the Riemannian
(Ricci flow) case, after replacing the assumption on the
nonnegativity of the bisectional curvature by the nonnegativity of
the curvature operator whenever necessary.

We first let $(M, g_{\abb}(x, t))$ ($m=\dim_\C M$) be a solution to
the K\"ahler-Ricci flow:
$$
\frac{\partial}{\partial t} g_{\abb}=-R_{\abb}.
$$
Let $\k_{\a\bbar}(x,t)$ be a Hermitian symmetric tensor defined on
$M\times [0,T]$, which is  deformed by the complex
Lichnerowicz-Laplacian heat equation (or L-heat equation in short):
$$
\left(\frac{\p}{\p t}-\D \right)\k_{\g\delbar}= R_{\beta
\bar{\a}\g\delbar}\k_{\a\bbar}-
\frac{1}{2}\left(R_{\g\bar{p}}k_{p\delbar}+
R_{p\delbar}\k_{\g\bar{p}}\right).
$$
Let $div(\k)_\a= g^{\g\delbar}\nabla_\g \k_{\a\delbar}$ and
$div(\k)_{\bbar}= g^{\g\delbar}\nabla_{\delbar}\k_{\g \bbar}$.
Consider the quantity
\begin{eqnarray*}
 Z& =&g^{\a\bbar}g^{\g\delbar}\left[\frac{1}{2} \lf
(\nabla_{\bbar}\nabla_\g +\nabla_\g \nabla_{\bbar}\ri)
\k_{\a\delbar}+R_{\a\delbar}\k_{\g\bbar}\right.\\
&\,&  +\left.\lf(\nabla_\g\k_{\a\delbar}V_{\bbar}+
\nabla_{\bbar}\k_{\a\delbar}V_{\g}\ri)+\k_{\a\delbar}V_{\bbar}V_\g
\ri]
 +\frac{K}{t}\\
&=&\frac12[g^{\a\bbar}\nabla_\bbar div(\k)_\a+g^{\g\delbar}\nabla_\g div(\k)_\delbar]\\
&\quad&+g^{\a\bbar}g^{\g\delbar}[R_{\a\delbar}\k_{\g\bbar}+\nabla_\g
\k_{\a\delbar}V_{\bbar}+
\nabla_{\bbar}\k_{\a\delbar}V_{\g}+\k_{\a\delbar}V_{\bbar}V_\g
]+\frac{K}{t}
\end{eqnarray*}
where $K$ is the trace of $\k_{\abb}$ with respect to
$g_{\a\bbar}(x,t)$.
 In \cite{NT} the following result,
which is the K\"ahler analogue of an earlier result in \cite{CH},
was showed by the maximum principle.

 \begin{theorem} \label{Theorem 4} Let $\k_{\abb}$ be a
Hermitian symmetric tensor satisfying the L-heat equation on
$M\times[0,T]$. Suppose $\k_{\abb}(x,0)\ge0$ (and satisfies some
growth assumptions in the case $M$ is noncompact). Then $Z\ge0$ on
$M\times(0,T]$ for any smooth vector field $V$ of type $(1,0)$.
\end{theorem}

The use of the maximum principle in the proof can be replaced by the
integration argument as in the proof of (\ref{4}). For any $T\ge
t_0>0$, in order to prove that $Z\ge 0$ at $t_0$ it suffices to show
that
 when $t=t_0$, $\int_M t^2Z h\, d\mu_t\ge 0$ for any
 compact-supported nonnegative function $h$. Now we solve the {\it
 conjugate heat equation} $\left(\frac{\partial }{\partial \tau} -\Delta
 +R\right)h(y,\tau)=0$ with $\tau=t_0-t$ and $h(y, \tau=0)=h(y)$,
 the given compact-supported function at $t_0$. By the perturbation
 argument we may as well as assume that $\k>0$. Then let
 $Z_m(y,t)=\inf_{V} Z(y, t)$. It was shown in \cite{NT} that
 $$ \lf(\frac{\p}{\p t}-\D \ri) Z_m =
Y_1+Y_2-2\frac{Z_m}{t}
$$
where \begin{eqnarray*}Y_1&=&\k_{\bar{p}q}\left(\D R_{p\bar{q}}
+R_{p\bar{q}\a\bbar}R_{\bar{\a}\beta} + \nabla_\a
R_{p\bar{q}}V_{\bar{\a}}+
\nabla_{\bar{\a}}R_{p\bar{q}}V_{\a}\right.\\
&\,& +\left. R_{p\bar{q}\a\bbar}V_{\bar{\a}}V_\beta +
\frac{R_{p\bar{q}}}{t}\right)
\end{eqnarray*}
and
\begin{eqnarray*}
Y_2&= &\k_{\g\bar{\a}}\lf[ \nabla_p V_{\bar{\g}}- R_{p\bar{\g}}
-\frac{1}{t}g_{p\bar{\g}}\ri]\lf[\nabla_{\bar{p}}V_\a- R_{\a\bar{p}}
-\frac{1}{t} g_{\bar{p}\a}\ri]\\
&\, &+
\k_{\g\bar{\a}}\nabla_{\bar{p}}V_{\bar{\g}}\nabla_pV_{\a}\\
&\ge& 0.
\end{eqnarray*}
Notice that in the above expressions, at every point $(y,t)$ the $V$
is the minimizing vector. This implies the monotonicity
$$
\frac{d}{d t}\int_M t^2Z_m h\, d\mu_t =t^2\int_M
\left(Y_1+Y_2\right)h\, d\mu_t\ge 0.
$$
Since $\lim_{t\to 0}t^2Z_m =0$, which is certainly the case if
$\Upsilon$ is smooth at $t=0$ and can be assumed so in general by
shifting $t$ with a $\e>0$, we have that $\int_M t^2Z_m h\, d\mu
|_{t=t_0}\ge 0$. This proof via the integration by parts implies the
following monotonicity formula.

\begin{proposition}\label{Proposition 5} Let $(M, g(t))$, $\k$ and $Z$ be as in Theorem
\ref{Theorem 4}. For any space-time point  $(x_0, t_0)$ with
$0<t_0\le T$, let $\ell(y, \tau)$ be the reduced distance function
with respect to $(x_0, t_0)$. Then
\begin{equation}\label{14}
\frac{d}{d t}\int_M t^2 Z_m \left(\frac{\exp(-\ell)}{(\pi
\tau)^m}\right)\, d\mu_t \ge t^2\int_M
\left(Y_1+Y_2\right)\left(\frac{\exp(-\ell)}{(\pi \tau)^m}\right)\,
d\mu_t\ge 0. \end{equation} In particular,
\begin{equation}\label{15}
t_0^2Z (x_0, t_0)\ge \int_0^{t_0}t^2\left(\int_M
\left(Y_1+Y_2\right)\left(\frac{\exp(-\ell)}{(\pi \tau)^m}\right)\,
d\mu_t\right)\, dt\ge 0. \end{equation}
\end{proposition}
Notice that (\ref{15}) sharpens the original Li-Yau-Hamilton
estimate of \cite{NT}, by encoding the rigidity (such as Hamilton's
result on eternal solutions), out of the equality case in the
Li-Yau-Hamilton estimate $Z\ge 0$, into the integral  of the right
hand side. The result holds for the Riemannian case if one uses
computation from \cite{CH}.

In \cite{N3}, the author discovered a new matrix Li-Yau-Hamilton
inequality for the K\"ahler-Ricci flow. (We also showed a family of
equations which connects this matrix inequality to Perelman's
entropy formula.) More precisely we showed that for any positive
solution $u$ to the {\it forward conjugate heat equation}
$\left(\frac{\p}{\p t}-\D -R\right)u=0$, we have that
\begin{equation}\label{16}
\k_{\abb} :=u\left(\nabla_\a\nabla_\bb \log
u+R_{\abb}+\frac{1}{t}g_{\abb}\right)\ge 0 \end{equation} under the
assumption that $(M, g(t))$ has bounded nonnegative bisectional
curvature.
 Using the above argument we can also
obtain a new monotonicity related to (\ref{16}). Indeed, tracing
(1.21) of \cite{N3} gives that
$$
\left(\frac{\p}{\p t}-\D\right)
Q=RQ-R_{\abb}\k_{\b\ba}-\frac{2}{t}Q+\frac{1}{u}|\k_{\abb}|^2+u\left|\nabla_\a\nabla_\b
\log u\right|^2+Y_3
$$
where $Q=g^{\abb}\k_{\abb}$ and
\begin{eqnarray*}
Y_3&=&u\left(\D R+|R_{\abb}|^2+\nabla_\a R\nabla_{\ba} \log
u+\nabla_\a\log u \nabla_{\ba} R\right.\\
&\, &+\left.R_{\abb}\nabla_\ba \log u\nabla_\b \log
u+\frac{1}{t}R\right)\\
&\ge& 0. \end{eqnarray*} Hence we have the following monotonicity
formula, noticing that
$$Y_4:=RQ-R_{\abb}\k_{\b\ba}\ge 0.$$

\begin{proposition}\label{Proposition 6} Let $(M, g(t))$ and $(x_0, t_0)$ be as in
Proposition \ref{Proposition 5}. Then
\begin{eqnarray}\label{17}
 &\,& \frac{d}{d t} \int_M t^2Q \left(\frac{\exp(-\ell)}{(\pi
\tau)^m}\right)\, d\mu_t \\
&\ge&
t^2\int_M\left(\frac{1}{u}|\k_{\abb}|^2+u\left|\nabla_\a\nabla_\b
\log u\right|^2+Y_3+Y_4\right)\left(\frac{\exp(-\ell)}{(\pi
\tau)^m}\right)\, d\mu_t\nonumber  \\
&\ge& 0. \nonumber
\end{eqnarray}
In particular,
\begin{eqnarray}\label{18}
 &\,&t_0^2 Q(x_0, t_0)\\
 &\ge &\int_0^{t_0} t^2\int_M\left(\frac{1}{u}|\k_{\abb}|^2+u\left|\nabla_\a\nabla_\b
\log u\right|^2+Y_3+Y_4\right)\left(\frac{\exp(-\ell)}{(\pi
\tau)^m}\right). \nonumber\end{eqnarray}
\end{proposition}
Again the advantage of the above monotonicity formula is that it
encodes the consequence on equality case (which is that $(M, g(t))$
is an gradient expanding soliton) into the  the right hand side
integral.

 Without Ricci flow, we can apply the similar argument to prove
Li-Yau's inequality and obtain a monotonicity formula. More
precisely, let $(M, g)$ ($n=\dim_\R M$) be a complete Riemannian
manifold with nonnegative Ricci curvature. Let $u(x,t)$ be a
positive solution to the heat equation on $M\times [0, T]$. Li and
Yau proved that
$$
\D \log u +\frac{n}{2t} \ge 0.
$$
Another way of  proving the above Li-Yau's inequality is through the
above integration by parts argument and the differential equation
$$
\left(\frac{\p}{\p t}-\D\right) Q =
\frac{2}{u}|\k_{ij}|^2-\frac{2}{t}Q +\frac{2}{u}R_{ij}\nabla_i
u\nabla_j u
$$
where
$$
\k_{ij}=\nabla_i\nabla_j u +\frac{u}{2t}g_{ij}-\frac{u_iu_j}{u}
$$
and $Q=g^{ij}\k_{ij}=u(\D \log u+\frac{n}{2t})$. This together with
Cheeger-Yau's theorem \cite{CY} on lower bound of the heat kernel,
gives the following monotonicity formula, which also give
characterization on the manifold if the equality holds somewhere for
some positive $u$.

\begin{proposition}\label{Proposition 7} Let $(M, g)$ be a complete Riemannian manifold
with non-negative Ricci curvature. Let $(x_0, t_0)$ be a space-time
point with $t_0>0$. Let $\tau=t_0-t$. Then
\begin{eqnarray}\label{19}
 &\,& \frac{d}{d t}\left(\int_M t^2 Q (y, t)
\hat H(x_0, y, \tau)\, d\mu(x)\right)\\
& \ge &2t^2\int_M \left(\left|\nabla_i\nabla_j \log
u+\frac{1}{2t}g_{ij}\right|^2+R_{ij}\nabla_i \log u \nabla_j \log
u\right)u \hat H\, d\mu\ge 0, \nonumber\end{eqnarray} where $\hat
H(x_0, y, \tau)= \frac{1}{(4\pi
\tau)^{\frac{n}{2}}}\exp\left(-\frac{d^2(x_0, y)}{4\tau}\right)$
with $d(x_0, y)$ being the distance function between $x_0$ and $y$.
In particular, we have that
\begin{eqnarray}\label{20}
 &\,& \left(u\D \log u+\frac{n}{2t}u\right)
(x_0, t_0)\\
&\ge &\frac{2}{t_0^2}\int_0^{t_0}t^2\int_M
\left(\left|\nabla_i\nabla_j \log
u+\frac{1}{2t}g_{ij}\right|^2+R_{ij}\nabla_i \log u \nabla_j \log
u\right)u \hat H.\nonumber
\end{eqnarray}
\end{proposition}

It is clear that (\ref{20}) improves the estimate of Li-Yau slightly
by providing the lower estimate, from which one can see easily that
the equality (for Li-Yau's estimate) holding  somewhere implies that
$M=\R^n$ (this was first observed in \cite{N1}, with the help of  an
entropy formula). The expression in the right hand side of
(\ref{17}) also appears in the linear entropy formula of \cite{N1}.

One can write down similar improving results for the Li-Yau type
estimate proved in \cite{N1}, which is a linear analogue of
Perelman's estimate $v_H\le 0$, and the one in \cite{N2}, which is a
linear version of Theorem 4 above. For example, when $M$ is a
complete Riemannian manifold with the nonnegative Ricci curvature,
if $u=H(x,y,t)=\frac{e^{-f}}{(4\pi t)^{\frac{n}{2}}}$, the
fundamental solution to the heat equation centered at $x$ at $t=0$,
letting $W=t(2\D f -|\nabla f|)+f-n$, we have that $W\le 0$. If
$\hat H$ is the `pseudo backward heat kernel' defined as in
Proposition \ref{Proposition 7} we have that
\begin{eqnarray*}
 &\,& \frac{d}{dt}\int_M (-W)u \hat H(x_0, y, \tau)\,
d\mu(y)\\
&=&2t\int_M\left(\left|\nabla_i \nabla_j
f-\frac{1}{2t}g_{ij}\right|^2+R_{ij}\nabla_i f\nabla_j f\right)u
\hat H(x_0, y, \tau)\, d\mu(y)\ge 0\end{eqnarray*} and $$
\left(-Wu\right)(x_0, t_0)\ge
2\int_0^{t_0}t\int_M\left(\left|\nabla_i \nabla_j
f-\frac{1}{2t}g_{ij}\right|^2+R_{ij}\nabla_i f\nabla_j f\right)u
\hat H.
$$
If we assume further that  $M$ is a complete K\"ahler manifold with
nonnegative bisectional curvature and $u(y,t)$ is a strictly
plurisubharmonic solution to the heat equation with $w=u_t$, then
$$
\frac{d}{ dt} \int_M t^2Z^w_m \hat H(x_0, y, t) \, d\mu(y)
=t^2\int_M Y_5 \hat H(x_0, y, t)\,\, d\mu(y)\ge 0,
$$
where
$$
Z^w_m(y, t)=\inf_{V\in T^{1, 0}M}\left(w_t+\nabla_\a w
V_{\ba}+\nabla_{\ba} w V_\a
+u_{\abb}V_{\ba}V_{\b}+\frac{w}{t}\right)
$$and
\begin{eqnarray*}
 Y_5 &=&u_{\gamma
\bar{\a}}\lf[\nabla_pV_{\bar{\gamma}}-\frac{1}{t}g_{p\bar{\gamma}}\ri]\lf[\nabla_{\bar{p}}V_{\a}
-\frac{1}{t}g_{\bar{p}\a}\ri]\\
&\, &+ u_{\gamma
\bar\a}\nabla_{\bar{p}}V_{\bar{\gamma}}\nabla_pV_{\a} +R_{\abb
s\bar{t}}u_{\bar{s}t}V_{\beta}V_{\bar{\a}}\\&\ge& 0
\end{eqnarray*}
with $V$ being the minimizing vector in the definition of  $ Z^w_m$.
In particular,
$$
\left(\frac{\p^2}{\p (\log t)^2} u(y, t) \right)(x_0,
t_0)\ge\int_0^{t_0}t^2\int_M Y_5 \hat H(x_0, y, t)\,\, d\mu(y)\, dt.
$$
This sharpens the logarithmic-convexity of \cite{N2}.

Finally we should remark that in all the  discussions above one can
replace the `pseudo backward heat kernel' $\hat H(y,t; x_0,
t_0)=\frac{\exp(-\frac{r^2(x_0, y)}{4(t_0-t)})}{(4\pi
(t_0-t))^{\frac{n}{2}}}$ (or $\frac{\exp(-\ell(y, \tau))}{(4\pi
\tau)^{\frac{n}{2}}}$, centered at $(x_0, t_0)$ in the case of Ricci
flow), which we wrote before as $\hat H(y, x_0, \tau)$ by  abusing
the notation, by the fundamental solution to the backward heat
equation (even by constant $1$ in the case of compact manifolds).
Also it still remains interesting on how to make effective uses of
these improved estimates, besides the rigidity results out of the
inequality being equality somewhere. There is also a small point
that should not be glossed over. When the manifold is complete
noncompact, one has to justify the validity of the Green's second
identity (for example in Proposition \ref{Proposition 7} we need to
justify that $\int_M \left(\hat H\D Q-Q\D \hat H\right)\, d\mu =0$).
This can be done when  $t_0$ is sufficiently small together with
integral estimates on the Li-Yau-Hamilton quantity (cf. \cite{CLN}).
The local monotonicity formula that shall be discussed in the next
section provides another way to avoid possible  technical
complications caused by the non-compactness.

\section{ Local monotonicity formulae}

In \cite{EKNT}, a very general scheme on localizing the monotonicity
formulae is developed. It is for any family of metrics evolved by
the equation $\frac{\p}{\p t} g_{ij}=-2\kappa_{ij}$. The
localization is through the so-called `heat ball'. More precisely
for a smooth positive space-time function $v$, which often is the
fundamental solution to the {\it backward conjugate heat equation}
or the `pseudo backward heat kernel' $\hat H(x_0, y,
\tau)=\frac{e^{-\frac{r^2(x_0, y)}{4\tau}}}{(4\pi
\tau)^{\frac{n}{2}}}$ (or $\frac{e^{-\ell(y, \tau)}}{(4\pi
\tau)^{\frac{n}{2}}}$ in the case of Ricci flow), with $\tau=t_0-t$,
one defines the `heat ball' by $E_r=\{(y, t)|\, v\ge r^{-n};  t<
t_0\}$. For all interesting cases we can check that $E_r$ is compact
for small $r$ (cf. \cite{EKNT}). Let $\psi_r=\log v +n\log r$. For
any `Li-Yau-Hamilton' quantity $\Q$ we define the local quantity:
$$
P(r):=\int_{E_r}\left(|\nabla \psi_r|^2+\psi_r({\text tr}_g
\kappa)\right)\Q\, d\mu_t \,dt.
$$
The finiteness of the integral can be verified via the localization
of Lemma \ref{Lemma 2}, a local gradient estimate. The general form
of the theorem, which is proved in Theorem 1 of \cite{EKNT}, reads
as the following.

\begin{theorem}\label{Theorem 8} Let $I(r)=\frac{P(r)}{r^n}$. Then
\begin{eqnarray}\label{21}
I(r_2)-I(r_1)&=&-\int_{r_1}^{r_2}
\frac{n}{r^{n+1}}\int_{E_r}\left[\left(\left(\frac{\p}{\p
t}+\D-{\text tr}_g
\kappa\right)v\right)\frac{\Q}{v}\right.\\
&\, &\left.+\psi_r\left(\frac{\p}{\p t}-\D \right)\Q\right] \,
d\mu_t \, dt\, dr. \nonumber\end{eqnarray}
\end{theorem}
It gives the monotonicity of $I(r)$ in the cases that $\Q \ge 0$,
which is ensured by the Li-Yau-Hamilton estimates in the case we
shall consider, and both $\left(\frac{\p}{\p t}+\D-{\text tr}_g
\kappa\right)v $ and $ \left(\frac{\p}{\p t}-\D \right)\Q$ are
nonnegative. The nonnegativity of  $\left(\frac{\p}{\p t}+\D-{\text
tr}_g \kappa\right)v $ comes for free if we chose $v$ to be the
`pseudo backward heat kernel'. The nonnegativity of $
\left(\frac{\p}{\p t}-\D \right)\Q$ follows from the computation,
which we may call as in \cite{N3} the {\it pre-Li-Yau-Hamilton
equation}, during the proof of the corresponding Li-Yau-Hamilton
estimate. Below we illustrate examples corresponding to the
monotonicity formulae derived in the previous section. These new
ones  expand the list of examples given in Section 4 of \cite{EKNT}.

 For the case of Ricci/K\"ahler-Ricci flow,  for
a fixed $(x_0, t_0)$, let $v=\frac{e^{-\ell(y, \tau)}}{(4\pi
\tau)^{\frac{n}{2}}}$, the `pseudo backward heat kernel', where
$\ell$ is the reduced distance centered at $(x_0, t_0)$.

 \begin{example}\label{Example 1} Let $Z_m$, $Y_1$ and $Y_2$ be as
in Proposition \ref{Proposition 5}. Let $\Q=t^2Z_m$. Then
$$
\frac{d}{d r}I(r)\le
-\frac{n}{r^{n+1}}\int_{E_r}\left[t^2\psi_r\left(Y_1+Y_2\right)\right]\,
d\mu_t \, dt\le 0
$$
and
$$
\Q(x_0, t_0)\ge
I(\bar{r})+\int_0^{\bar{r}}\frac{n}{r^{n+1}}\int_{E_r}\left[t^2\psi_r\left(Y_1+Y_2\right)\right]\,
d\mu_t \, dt\, dr.
$$
\end{example}

\begin{example}\label{Example 2} Let $u$, $\Q=t^2Q$, $\Upsilon_{\abb}$, $Y_3$ and
$Y_4$ be as in Proposition \ref{Proposition 6}. Then
$$
\frac{d}{d r}I(r)\le
-\frac{n}{r^{n+1}}\int_{E_r}t^2\psi_r\left(\frac{1}{u}|\k_{\abb}|^2+u\left|\nabla_\a\nabla_\b
\log u\right|^2+Y_3+Y_4\right)\le 0
$$and
$$ \Q(x_0, t_0)\ge I(\bar{r})+\int_0^{\bar{r}}\frac{n}{r^{n+1}}\int_{E_r}t^2\psi_r\left(\frac{1}{u}|\k_{\abb}|^2+u\left|\nabla_\a\nabla_\b
\log u\right|^2+Y_3+Y_4\right). $$
\end{example}

For the fixed metric case, we may choose either $v=H(x_0, y, \tau)$,
the {\it backward heat kernel} or $v=\hat H(x_0, y,
\tau)=\frac{e^{-\frac{d^2(x_0, y)}{4\tau}}}{(4\pi
\tau)^{\frac{n}{2}}}$, the `pseudo backward heat kernel'.

\begin{example}\label{Example 3} Let $u$ and  $Q$ be as in Proposition \ref{Proposition 7}. Let
$\Q=t^2 Q$ and $f=\log u$. Then
$$
\frac{d}{d r}I(r)\le -\frac{2n}{r^{n+1}}\int_{E_r}t^2u\psi_r
\left(\left|\nabla_i\nabla_j
f+\frac{1}{2t}g_{ij}\right|^2+R_{ij}\nabla_i f \nabla_j f\right)\,
d\mu\, dt\le 0
$$
and $$ \Q(x_0, t_0)\ge I(\bar{r})+\int_0^{\bar{r}}
\frac{2n}{r^{n+1}}\int_{E_r}t^2u\psi_r\left(\left|\nabla_i\nabla_j
f+\frac{1}{2t}g_{ij}\right|^2+R_{ij}\nabla_i f \nabla_j f\right) .$$
\end{example}

\begin{example}\label{Example 4} Let $u=\frac{e^{-f}}{(4\pi t)^{\frac{n}{2}}}$ be
the fundamental solution to the (regular) heat equation. Let
$W=t(2\D f-|\nabla f|^2)+f-n$ and $\Q=-uW$. Then
$$
\frac{d}{d r}I(r)\le
-\frac{2n}{r^{n+1}}\int_{E_r}tu\psi_r\left(\left|\nabla_i \nabla_j
f-\frac{1}{2t}g_{ij}\right|^2+R_{ij}\nabla_i f\nabla_j f\right) \,
d\mu\, dt\le 0
$$
and $$
 \Q(x_0, t_0)\ge I(\bar{r})+\int_0^{\bar{r}}\frac{2n}{r^{n+1}}\int_{E_r}tu\psi_r\left(\left|\nabla_i \nabla_j
f-\frac{1}{2t}g_{ij}\right|^2+R_{ij}\nabla_i f\nabla_j f\right).
$$Note that this provides another localization of entropy other than
the one in \cite{N3} (see also \cite{CLN}).
\end{example}

\begin{example}
\label{Example 5} Let $M$ be a complete K\"ahler manifold with
nonnegative bisectional curvature. Let $u$, $Z^w_m$ and $Y_5$ be as
in the last case considered in Section 4. Let $\Q=t^2Z^w_m$. Then
$$
\frac{d}{d r}I(r)\le -\frac{n}{r^{n+1}}\int_{E_r} t^2 Y_5 \psi_r \,
d\mu \, dt
$$ and
$$
\left(\frac{\p^2}{\p (\log t)^2} u(x, t) \right)(x_0, t_0)\ge
I(\bar{r})+\int_0^{\bar{r}}\frac{n}{r^{n+1}}\int_{E_r} t^2 Y_5
\psi_r \, d\mu \, dt\, dr.
$$
\end{example}

\medskip

{\it Acknowledgement}. We would like to thank  Ben Chow and Peng Lu
for continuously pressing us on a understandable proof of (\ref{4}).
We started to seriously  work on it after the  visit to  Klaus Ecker
in August and a stimulating discussion with him. We would like to
thank him for that, as well as Dan Knopf and Peter Topping for
discussions on a related issue.

\bibliographystyle{amsalpha}

\end{document}